\documentclass[thmsa,12pt]{article}
\usepackage{amsmath}
\usepackage{amssymb}
\usepackage{sw20lart}

\input{tcilatex}

\begin{document}

\title{Extensors in Geometric Algebra}
\author{{\footnotesize A. M. Moya}$^{2},${\footnotesize V. V. Fern\'{a}ndez}$%
^{1}${\footnotesize \ and W. A. Rodrigues Jr}.$^{1}${\footnotesize \ } \\
$^{1}\hspace{-0.1cm}${\footnotesize Institute of Mathematics, Statistics and
Scientific Computation}\\
{\footnotesize \ IMECC-UNICAMP CP 6065}\\
{\footnotesize \ 13083-859 Campinas, SP, Brazil }\\
{\footnotesize e-mail: walrod@ime.unicamp.br}\\
{\footnotesize \ }$^{2}${\footnotesize Department of Mathematics, University
of Antofagasta, Antofagasta, Chile} \\
{\footnotesize e-mail: mmoya@uantof.cl }}
\maketitle

\begin{abstract}
This paper, the third in a series of eight introduces some of the basic
concepts of the theory of \emph{extensors} needed for our formulation of the
differential geometry of smooth manifolds . Key notions such as the
extension and generalization operators of a given linear operator (a ($1,1$%
)-extensor) acting on a real vector space $V$ are introduced and studied in
details. Also, we introduce the notion of the determinant of a ($1,1$%
)-extensor and the concepts of standard and metric Hodge (star) operators,
disclosing a non trivial and useful relation between them.
\end{abstract}

\tableofcontents

\section{Introduction}

In this paper, the third in a series of eight dealing formulation of the
differential geometry of smooth manifolds with Clifford (geometric) algebra
methods we recall some basic notions of the theory of extensors thus
completing the presentation of the algebraic notions necessary for the
remaining papers of the series. Here, in Section 2 we introduce $k$%
-extensors and $(p,q)$-extensors. Section 3 deals with projector operators.
Section 4 introduces the key concept of the extension operator (or exterior
power operator) of a \textit{given} linear operator $t$ (i.e., a $(1,1)$%
-extensor) on a real vector space $V$. In section 5 the standard and metric
adjoint operators are recalled as extensor operators. Section 6 introduces
the generalization operator of $t$, an important concept in our formulation
of the differential geometry of arbitrary (smooth) manifolds which is
presented in sequel papers of this series. \ \ Section 7 introduces the
concept of the determinant of a $(1,1)$-extensor and Section 8 shows through
some applications the concept of extensors in action. There, the standard
and metric Hodge operators are introduced and an explicit formula connecting
them is derived. Moreover, another important formula is derived which
relates the metric Hodge operators corresponding to two distorted metric
structures on $V$. In Section 9 \ we present our conclusions.

\section{General $k$-Extensors}

Let $\bigwedge_{1}^{\diamond }V,\ldots ,$ $\bigwedge_{k}^{\diamond }V$ be $k$
subspaces of $\bigwedge V$ such that each of them is \emph{any} sum of
homogeneous subspaces of $\bigwedge V,$ and $\bigwedge^{\diamond }V$ is
either any sum of homogeneous subspaces of $\bigwedge V$ or the trivial
subspace $\{0\}.$ A multilinear mapping from the cartesian product $%
\bigwedge_{1}^{\diamond }V\times \cdots \times \bigwedge_{k}^{\diamond }V$
to $\bigwedge^{\diamond }V$ will be called a general $k$-extensor over $V,$
i.e., $t:\bigwedge_{1}^{\diamond }V\times \cdots \times
\bigwedge_{k}^{\diamond }V\rightarrow \bigwedge^{\diamond }V$ such that for
any $\alpha _{j},\alpha _{j}^{\prime }\in \mathbb{R}$ and $%
X_{j},X_{j}^{\prime }\in \bigwedge_{j}^{\diamond }V,$ 
\begin{equation}
t(\ldots ,\alpha _{j}X_{j}+\alpha _{j}^{\prime }X_{j}^{\prime },\ldots
)=\alpha _{j}t(\ldots ,X_{j},\ldots )+\alpha _{j}^{\prime }t(\ldots
,X_{j}^{\prime },\ldots ),  \label{GE.1}
\end{equation}
for each $j$ with $1\leq j\leq k.$

It should be noticed that the linear operators on $V,$ $\bigwedge^{p}V$ or $%
\bigwedge V$ which appear in ordinary linear algebra are particular cases of 
$1$-extensors over $V.$ Note also that a covariant $k$-tensor over $V$ is
just a $k$-extensor over $V.$ On this way, the concept of general $k$%
-extensor generalizes and unifies both of the concepts of linear operator
and of covariant $k$-tensor. These mathematical objects are of the same
nature!

The set of general $k$-extensors over $V,$ denoted by $k$-$%
ext(\bigwedge_{1}^{\diamond }V,\ldots ,\bigwedge_{k}^{\diamond
}V;\bigwedge^{\diamond }V),$ has a natural structure of real vector space.
Its dimension is clearly given by 
\begin{equation}
\dim k\text{-}ext(\bigwedge_{1}^{\diamond }V,\ldots ,\bigwedge_{k}^{\diamond
}V;\bigwedge^{\diamond }V)=\dim \bigwedge_{1}^{\diamond }V\cdots \dim
\bigwedge_{k}^{\diamond }V\dim \bigwedge^{\diamond }V.  \label{GE.2}
\end{equation}%
We shall need to consider only some particular cases of these general $k$%
-extensors over $V.$ So, special names and notations will be given for them.

We will equip $V$ with an arbitrary (but fixed once and for all) euclidean
metric $G_{E}$, and denote the scalar product of multivectors $X,Y\in
\bigwedge V$ with respect to the euclidean metric structure $(V,G_{E}),$ $%
X\cdot Y$ instead of the more detailed notation $X\underset{G_{E}}{\cdot }Y$.

Let $\{e_{j}\}$ be any basis for $V,$ and $\{e^{j}\}$ be its euclidean
reciprocal basis for $V,$ i.e., $e_{j}\cdot e^{k}=\delta _{j}^{k}$.

\subsection{$(p,q)$-Extensors}

Let $p$ and $q$ be two integer numbers with $0\leq p,q\leq n$. A linear
mapping which sends $p$-vectors to $q$-vectors will be called a $(p,q)$-%
\emph{extensor} over $V.$ The space of these objects, namely $1$-$%
ext(\bigwedge^{p}V;\bigwedge^{q}V),$ will be denoted by $ext_{p}^{q}(V)$ for
short. By using Eq.(\ref{GE.2}) we get 
\begin{equation}
\dim ext_{p}^{q}(V)=\binom{n}{p}\binom{n}{q}.  \label{GE.3}
\end{equation}

For instance, we see that the $(1,1)$-extensors over $V$ are just the
well-known linear operators on $V.$

The $\binom{n}{p}\binom{n}{q}$ extensors belonging to $ext_{p}^{q}(V),$
namely $\varepsilon ^{j_{1}\ldots j_{p};k_{1}\ldots k_{q}},$ defined by 
\begin{equation}
\varepsilon ^{j_{1}\ldots j_{p};k_{1}\ldots k_{q}}(X)=(e^{j_{1}}\wedge
\ldots \wedge e^{j_{p}})\cdot Xe^{k_{1}}\wedge \ldots \wedge e^{k_{q}}
\label{GE.3a}
\end{equation}
define a $(p,q)$-extensor basis for $ext_{p}^{q}(V).$

Indeed, the extensors given by Eq.(\ref{GE.3a}) are linearly independent,
and for each $t\in ext_{p}^{q}(V)$ there exist $\binom{n}{p}\binom{n}{q}$
real numbers, say $t_{j_{1}\ldots j_{p};k_{1}\ldots k_{q}},$ given by 
\begin{equation}
t_{j_{1}\ldots j_{p};k_{1}\ldots k_{q}}=t(e_{j_{1}}\wedge \ldots \wedge
e_{j_{p}})\cdot (e_{k_{1}}\wedge \ldots \wedge e_{k_{q}})  \label{GE.3b}
\end{equation}%
such that 
\begin{equation}
t=\frac{1}{p!}\frac{1}{q!}t_{j_{1}\ldots j_{p};k_{1}\ldots k_{q}}\varepsilon
^{j_{1}\ldots j_{p};k_{1}\ldots k_{q}}.  \label{GE.3c}
\end{equation}%
Such $t_{j_{1}\ldots j_{p};k_{1}\ldots k_{q}}$ will be called the $%
j_{1}\ldots j_{p};k_{1}\ldots k_{q}$-th \emph{covariant components} of $t$
with respect to the $(p,q)$-\emph{extensor basis} $\{\varepsilon
^{j_{1}\ldots j_{p};k_{1}\ldots k_{q}}\}.$

Of course, there are still other kinds of $(p,q)$-extensor bases for $%
ext_{p}^{q}(V)$ besides the one given by Eq.(\ref{GE.3a}) which can be
constructed from the vector bases $\{e_{j}\}$ and $\{e^{j}\}.$ The total
number of these different kinds of $(p,q)$-extensor bases for $%
ext_{p}^{q}(V) $ are $2^{p+q}$.

Now, if we take the basis $(p,q)$-extensors $\varepsilon _{j_{1}\ldots
j_{p};k_{1}\ldots k_{q}}$ and the real numbers $t^{j_{1}\ldots
j_{p};k_{1}\ldots k_{q}}$ defined by 
\begin{eqnarray}
\varepsilon _{j_{1}\ldots j_{p};k_{1}\ldots k_{q}}(X) &=&(e_{j_{1}}\wedge
\ldots \wedge e_{j_{p}})\cdot Xe_{k_{1}}\wedge \ldots \wedge e_{k_{q}},
\label{GE.3d} \\
t^{j_{1}\ldots j_{p};k_{1}\ldots k_{q}} &=&t(e^{j_{1}}\wedge \ldots \wedge
e^{j_{p}})\cdot (e^{k_{1}}\wedge \ldots \wedge e^{k_{q}}),  \label{GE.3e}
\end{eqnarray}%
we get an expansion formula for $t\in ext_{p}^{q}(V)$ analogous to that
given by Eq.(\ref{GE.3c}), i.e., 
\begin{equation}
t=\frac{1}{p!}\frac{1}{q!}t^{j_{1}\ldots j_{p};k_{1}\ldots k_{q}}\varepsilon
_{j_{1}\ldots j_{p};k_{1}\ldots k_{q}}.  \label{GE.3f}
\end{equation}%
Such $t^{j_{1}\ldots j_{p};k_{1}\ldots k_{q}}$ are called the $j_{1}\ldots
j_{p};k_{1}\ldots k_{q}$-th \emph{contravariant components} of $t$ with
respect to the $(p,q)$-\emph{extensor basis} $\{\varepsilon _{j_{1}\ldots
j_{p};k_{1}\ldots k_{q}}\}.$

\subsection{ Extensors}

A linear mapping which sends multivectors to multivectors will be simply
called an \emph{extensor} over $V.$ They are the linear operators on $%
\bigwedge V.$ For the space of extensors over $V,$ namely $1$-$ext(\bigwedge
V;\bigwedge V),$ we will use the short notation $ext(V).$ By using Eq.(\ref%
{GE.2}) we get 
\begin{equation}
\dim ext(V)=2^{n}2^{n}.  \label{GE.4}
\end{equation}

For instance, we will see that the so-called Hodge star operator is just a
well-defined extensor over $V$ which can be thought as an exterior direct
sum of $(p,n-p)$-extensor over $V.$ The extended (or exterior power) of $%
t\in ext_{1}^{1}(V)$ is just an extensor over $V,$ i.e., $\underline{t}\in
ext(V).$

There are $2^{n}2^{n}$ extensors over $V,$ namely $\varepsilon ^{J;K}$,
given by\footnote{$J$ and $K$ are colective indices. Recall, for example,
that: $e_{J}=1,e_{j_{1}},e_{j_{1}}\wedge e_{j_{2}},\ldots $($%
e^{J}=1,e^{j_{1}},e^{j_{1}}\wedge e^{j_{2}},\ldots $) and $\nu
(J)=0,1,2,\ldots $ for $J=\emptyset ,j_{1},j_{1}j_{2},\ldots ,$ where all
index $j_{1},j_{2},\ldots $ runs from $1$ to $n.$} 
\begin{equation}
\varepsilon ^{J;K}(X)=(e^{J}\cdot X)e^{K}  \label{GE.4a}
\end{equation}%
which can be used to introduce an \emph{\ extensor} \emph{bases }for $%
ext(V). $

In fact they are linearly independent, and for each $t\in ext(V)$ there
exist $2^{n}2^{n}$ real numbers, say $t_{J;K},$ given by 
\begin{equation}
t_{J;K}=t(e_{J})\cdot e_{K}  \label{GE.4b}
\end{equation}
such that 
\begin{equation}
t=\underset{J}{\sum }\underset{K}{\sum }\frac{1}{\nu (J)!}\frac{1}{\nu (K)!}%
t_{J;K}\varepsilon ^{J;K}.  \label{GE.4c}
\end{equation}
Such $t_{J;K}$ will be called the $J;K$-th\emph{\ covariant components} of $%
t $ with respect to the \emph{extensor bases} $\{\varepsilon ^{J;K}\}.$

We notice that exactly $(2^{n+1}-1)^{2}$ extensor bases for $ext(V)$ can be
constructed from the basis vectors $\{e_{j}\}$ and $\{e^{j}\}.$ For
instance, whenever the basis extensors $\varepsilon _{J;K}$ and the real
numbers $t^{J;K}$ defined by 
\begin{eqnarray}
\varepsilon _{J;K}(X) &=&(e_{J}\cdot X)e_{K},  \label{GE.4d} \\
t^{J;K} &=&t(e^{J})\cdot e^{K}  \label{GE.4e}
\end{eqnarray}%
are used, an expansion formula for $t\in ext(V)$ analogous to that given by
Eq.(\ref{GE.4c}) can be obtained, i.e., 
\begin{equation}
t=\underset{J}{\sum }\underset{K}{\sum }\frac{1}{\nu (J)!}\frac{1}{\nu (K)!}%
t^{J;K}\varepsilon _{J;K}.  \label{GE.4f}
\end{equation}%
Such $t^{J;K}$ are called the $J;K$-th \emph{contravariant components} of $t$
with respect to the \emph{extensor bases} $\{\varepsilon _{J;K}\}$.

\subsection{Elementary $k$-Extensors}

A multilinear mapping which takes $k$-uple of vectors into $q$-vectors will
be called an \emph{elementary }$k$-\emph{extensor} over $V$ of degree $q.$
The space of these objects, namely $k$-$ext(V,\ldots ,V;\bigwedge^{q}V),$
will be denoted by $k$-$ext^{q}(V).$ It is easy to verify (using Eq.(\ref%
{GE.2})) that 
\begin{equation}
\dim k\text{-}ext^{q}(V)=n^{k}\binom{n}{q}.  \label{GE.5}
\end{equation}%
It should be noticed that an\emph{\ }elementary $k$-extensor over $V$ of
degree $0$ is just a \emph{covariant }$k$-\emph{tensor} over $V,$ i.e., $k$-$%
ext^{0}(V)\equiv T_{k}(V).$ It is easily realized that $1$-$ext^{q}(V)\equiv
ext_{1}^{q}(V).$

The elementary $k$-extensors of degrees $0,1,2,\ldots $ etc. are sometimes
said to be \emph{scalar}, \emph{vector,} \emph{bivector}, $\ldots $ etc. 
\emph{elementary }$k$-\emph{extensors.}

The $n^{k}\binom{n}{q}$ elementary $k$-extensors of degree $q$ belonging to $%
k$-$ext^{q}(V),$ namely $\varepsilon ^{j_{1},\ldots ,j_{k};k_{1}\ldots
k_{q}},$ given by 
\begin{equation}
\varepsilon ^{j_{1},\ldots ,j_{k};k_{1}\ldots k_{q}}(v_{1},\ldots
,v_{k})=(v_{1}\cdot e^{j_{1}})\ldots (v_{k}\cdot e^{j_{k}})e^{k_{1}}\wedge
\ldots \wedge e^{k_{q}}  \label{GE.5a}
\end{equation}
define elementary \emph{basis vectors}, ( i.e., $k$-extensor of degree $q$)
for $k$-$ext^{q}(V)$.

In fact they are linearly independent, and for all $t\in k$-$ext^{q}(V)$
there are $n^{k}\binom{n}{q}$ real numbers, say $t_{j_{1},\ldots
,j_{k};k_{1}\ldots k_{q}},$ given by 
\begin{equation}
t_{j_{1},\ldots ,j_{k};k_{1}\ldots k_{q}}=t(e_{j_{1}},\ldots
,e_{j_{k}})\cdot (e_{k_{1}}\wedge \ldots \wedge e_{k_{q}})  \label{GE.5b}
\end{equation}
such that 
\begin{equation}
t=\frac{1}{q!}t_{j_{1},\ldots ,j_{k};k_{1}\ldots k_{q}}\varepsilon
^{j_{1},\ldots ,j_{k};k_{1}\ldots k_{q}}.  \label{GE.5c}
\end{equation}
Such $t_{j_{1},\ldots ,j_{k};k_{1}\ldots k_{q}}$ will be called the $%
j_{1},\ldots ,j_{k};k_{1}\ldots k_{q}$-th \emph{covariant components }of $t$
with respect to the\emph{\ bases} $\{\varepsilon ^{j_{1},\ldots
,j_{k};k_{1}\ldots k_{q}}\}.$

We notice that exactly $2^{k+q}$ elementary $k$-extensors of degree $q$
bases for $k$-$ext^{q}(V)$ can be constructed from the vector bases $%
\{e_{j}\}$ and $\{e^{j}].$ For instance, we may define $\varepsilon
_{j_{1},\ldots ,j_{k};k_{1}\ldots k_{q}}$ (the basis elementary $k$-extensor
of degree $q$) and the real numbers $t^{j_{1},\ldots ,j_{k};k_{1}\ldots
k_{q}}$ by 
\begin{eqnarray}
\varepsilon _{j_{1},\ldots ,j_{k};k_{1}\ldots k_{q}}(v_{1},\ldots ,v_{k})
&=&(v_{1}\cdot e_{j_{1}})\ldots (v_{k}\cdot e_{j_{k}})e_{k_{1}}\wedge \ldots
\wedge e_{k_{q}},  \label{GE.5d} \\
t^{j_{1},\ldots ,j_{k};k_{1}\ldots k_{q}} &=&t(e^{j_{1}},\ldots
,e^{j_{k}})\cdot (e^{k_{1}}\wedge \ldots \wedge e^{k_{q}}).  \label{GE.5e}
\end{eqnarray}%
Then, we also have other expansion formulas for $t\in k$-$ext^{q}(V)$
besides that given by Eq.(\ref{GE.5c}), e.g., 
\begin{equation}
t=\frac{1}{q!}t^{j_{1},\ldots ,j_{k};k_{1}\ldots k_{q}}\varepsilon
_{j_{1},\ldots ,j_{k};k_{1}\ldots k_{q}}.  \label{GE.5f}
\end{equation}%
Such $t^{j_{1},\ldots ,j_{k};k_{1}\ldots k_{q}}$ are called the $%
j_{1},\ldots ,j_{k};k_{1}\ldots k_{q}$-th \emph{contravariant components} of 
$t$ with respect to the \emph{bases} $\{\varepsilon _{j_{1},\ldots
,j_{k};k_{1}\ldots k_{q}}\}.$

\section{Projectors}

Let $\bigwedge^{\diamond }V$ be either any sum of homogeneous subspaces%
\footnote{%
Note that for such a subspace $\bigwedge^{\diamond }V$ there are $\nu $
integers $p_{1,}\ldots ,p_{\nu }$ ($0\leq p_{1}<\cdots <p_{\nu }\leq n$)
such that $\bigwedge_{1}^{\diamond }V=\bigwedge^{p_{1}}V+\cdots
+\bigwedge^{p_{\nu }}V.$} of $\bigwedge V$ or the trivial subspace $\{0\}.$
Associated to $\bigwedge^{\diamond }V,$ a noticeable extensor from $%
\bigwedge V$ to $\bigwedge^{\diamond }V,$ namely $\left\langle \left.
{}\right. \right\rangle _{\bigwedge^{\diamond }V},$ can defined by 
\begin{equation}
\left\langle X\right\rangle _{\Lambda ^{\diamond }V}=\left\{ 
\begin{array}{cc}
\left\langle X\right\rangle _{p_{1}}+\cdots +\left\langle X\right\rangle
_{p\nu }, & \text{if }\bigwedge^{\diamond }V=\bigwedge^{p_{1}}V+\cdots
+\bigwedge^{p_{\nu }}V \\ 
0, & \text{if }\bigwedge^{\diamond }V=\{0\}%
\end{array}
\right. .  \label{P.1}
\end{equation}
Such $\left\langle \left. {}\right. \right\rangle _{\bigwedge^{\diamond
}V}\in 1$-$ext(\bigwedge V;\bigwedge^{\diamond }V)$ will be called the $%
\bigwedge^{\diamond }V$-\emph{projector extensor}.

We notice that if $\bigwedge^{\diamond }V$ is any homogeneous subspace of $%
\bigwedge V,$ i.e., $\bigwedge^{\diamond }V=\bigwedge^{p}V,$ then the
projector extensor is reduced to the so-called $p$-\emph{part operator},
i.e., $\left\langle \left. {}\right. \right\rangle _{\bigwedge^{\diamond
}V}=\left\langle \left. {}\right. \right\rangle _{p}.$

We now summarize the fundamental properties for the $\bigwedge^{\diamond }V$%
-projector extensors.

Let $\bigwedge_{1}^{\diamond }V$ and $\bigwedge_{2}^{\diamond }V$ be two
subspaces of $\bigwedge V.$ If each of them is either any sum of homogeneous
subspaces of $\bigwedge V$ or the trivial subspace $\{0\},$ then 
\begin{eqnarray}
\left\langle \left\langle X\right\rangle _{\bigwedge_{1}^{\diamond
}V}\right\rangle _{\bigwedge_{2}^{\diamond }V} &=&\left\langle
X\right\rangle _{\bigwedge_{1}^{\diamond }V\cap \bigwedge_{2}^{\diamond }V}
\label{P.2a} \\
\left\langle X\right\rangle _{\bigwedge_{1}^{\diamond }V}+\left\langle
X\right\rangle _{\bigwedge_{2}^{\diamond }V} &=&\left\langle X\right\rangle
_{\bigwedge_{1}^{\diamond }V\cup \bigwedge_{2}^{\diamond }V}.  \label{P.2b}
\end{eqnarray}

Let $\bigwedge^{\diamond }V$ be either any sum of homogeneous subspaces of $%
\bigwedge V$ or the trivial subspace $\{0\}.$ Then, it holds 
\begin{equation}
\left\langle X\right\rangle _{\bigwedge^{\diamond }V}\cdot Y=X\cdot
\left\langle Y\right\rangle _{\bigwedge^{\diamond }V}..  \label{P.2c}
\end{equation}

We see that the concept of $\bigwedge^{\diamond }V$-projector extensor is
just a natural generalization of the concept of $p$-part operator.

\section{The Extension Operator}

Let $\{e_{j}\}$ be any basis for $V,$ and $\{\varepsilon ^{j}\}$ be its dual
basis for $V^{*}.$ As we know, $\{\varepsilon ^{j}\}$ is the unique $1$-form
basis associated to the vector basis $\{e_{j}\}$ such that $\varepsilon
^{j}(e_{i})=\delta _{i}^{j}.$ The linear mapping $ext_{1}^{1}(V)\ni t\mapsto 
\underline{t}\in ext(V)$ such that for any $X\in \bigwedge V$ and $X=X_{0}+%
\overset{n}{\underset{k=1}{\dsum }}X_{k},$ then 
\begin{equation}
\underline{t}(X)=X_{0}+\overset{n}{\underset{k=1}{\dsum }}\frac{1}{k!}%
X_{k}(\varepsilon ^{j_{1}},\ldots ,\varepsilon ^{j_{k}})t(e_{j_{1}})\wedge
\ldots \wedge t(e_{j_{k}})  \label{EO.1}
\end{equation}
will be called the \emph{extension operator}. We call $\underline{t}$ the 
\emph{extended} of $t.$ It is the well-known outermorphism of $t$ in
ordinary linear algebra.

The extension operator is well-defined since it does not depend on the
choice of $\{e_{j}\}.$

We summarize now the basic properties satisfied by the extension operator.

\textbf{e1} The extension operator is grade-preserving, i.e., 
\begin{equation}
\text{if }X\in \bigwedge^{p}V,\text{ then }\underline{t}(X)\in
\bigwedge^{p}V.  \label{EO.1a}
\end{equation}%
It is an obvious result which follows from Eq.(\ref{EO.1}).

\textbf{e2} For any $\alpha \in \mathbb{R},$ $v\in V$ and $v_{1}\wedge
\ldots \wedge v_{k}\in \bigwedge^{k}V$, 
\begin{eqnarray}
\underline{t}(\alpha ) &=&\alpha ,  \label{EO.2a} \\
\underline{t}(v) &=&t(v),  \label{EO.2b} \\
\underline{t}(v_{1}\wedge \ldots \wedge v_{k}) &=&t(v_{1})\wedge \ldots
\wedge t(v_{k}).  \label{EO.2c}
\end{eqnarray}

\textbf{Proof}

The first statement trivially follows from Eq.(\ref{EO.1}). The second ones
can easily be deduced from Eq.(\ref{EO.1}) by recalling the elementary
expansion formula for vectors and the linearity of extensors. In order to
prove the third statement we use the formulas : $v_{1}\wedge \ldots \wedge
v_{k}(\omega ^{1},\ldots ,\omega ^{k})=\epsilon ^{i_{1}\ldots i_{k}}\omega
^{1}(v_{i_{1}})\ldots \omega ^{k}(v_{i_{k}})$ and $w_{i_{1}}\wedge \ldots
\wedge w_{i_{k}}=\epsilon _{i_{1}\ldots i_{k}}w_{1}\wedge \ldots \wedge
w_{k},$ where $v_{1},\ldots ,v_{k}\in V,$ $w_{1},\ldots ,w_{k}\in V$ and $%
\omega ^{1},\ldots ,\omega ^{k}\in V^{\ast },$ and the combinatorial formula 
$\epsilon ^{i_{1}\ldots i_{k}}\epsilon _{i_{1}\ldots i_{k}}=k!.$ From Eq.(%
\ref{EO.1}) by recalling the elementary expansion formula for vectors and
the linearity of extensors we have that 
\begin{eqnarray*}
\underline{t}(v_{1}\wedge \ldots \wedge v_{k}) &=&\frac{1}{k!}v_{1}\wedge
\ldots \wedge v_{k}(\varepsilon ^{j_{1}},\ldots ,\varepsilon
^{j_{k}})t(e_{j_{1}})\wedge \ldots \wedge t(e_{j_{k}}) \\
&=&\frac{1}{k!}\epsilon ^{i_{1}\ldots i_{k}}\varepsilon
^{j_{1}}(v_{i_{1}})\ldots \varepsilon ^{j_{k}}(v_{i_{k}})t(e_{j_{1}})\wedge
\ldots \wedge t(e_{j_{k}}) \\
&=&\frac{1}{k!}\epsilon ^{i_{1}\ldots i_{k}}t(v_{i_{1}})\wedge \ldots \wedge
t(v_{i_{k}}), \\
&=&t(v_{1})\wedge \ldots \wedge t(v_{k}).\blacksquare
\end{eqnarray*}

\textbf{e3} For any $X,Y\in \bigwedge V$, 
\begin{equation}
\underline{t}(X\wedge Y)=\underline{t}(X)\wedge \underline{t}(Y).
\label{EO.3}
\end{equation}%
Eq.(\ref{EO.3}) an immediate result which follows from Eq.(\ref{EO.2c}).

We emphasize that the three fundamental properties as given by Eq.(\ref%
{EO.2a}), Eq.(\ref{EO.2b}) and Eq.(\ref{EO.3}) together are completely
equivalent to the \emph{extension procedure} as defined by Eq.(\ref{EO.1}).

We present next some important properties for the extension operator.

\textbf{e4} Let us take $s,t\in ext_{1}^{1}(V).$ Then, the following result
holds 
\begin{equation}
\underline{s\circ t}=\underline{s}\circ \underline{t}.  \label{EO.4}
\end{equation}

\textbf{Proof}

It is enough to present the proofs for scalars and simple $k$-vectors.

For $\alpha \in \mathbb{R},$ by using Eq.(\ref{EO.2a}) we get 
\begin{equation*}
\underline{s\circ t}(\alpha )=\alpha =\underline{s}(\alpha )=\underline{s}(%
\underline{t}(\alpha ))=\underline{s}\circ \underline{t}(\alpha ).
\end{equation*}

For a simple $k$-vector $v_{1}\wedge \ldots \wedge v_{k}\in \bigwedge^{k}V,$
by using Eq.(\ref{EO.2c}) we get 
\begin{eqnarray*}
\underline{s\circ t}(v_{1}\wedge \ldots \wedge v_{k}) &=&s\circ
t(v_{1})\wedge \ldots \wedge s\circ t(v_{k})=s(t(v_{1}))\wedge \ldots \wedge
s(t(v_{k})) \\
&=&\underline{s}(t(v_{1})\wedge \ldots \wedge t(v_{k}))=\underline{s}(%
\underline{t}(v_{1}\wedge \ldots \wedge v_{k})), \\
&=&\underline{s}\circ \underline{t}(v_{1}\wedge \ldots \wedge v_{k}).
\end{eqnarray*}

Next, we can easily generalize to multivectors due to the linearity of
extensors. It yields 
\begin{equation*}
\underline{s\circ t}(X)=\underline{s}\circ \underline{t}(X).\blacksquare
\end{equation*}

\textbf{e5} Let us take $t\in ext_{1}^{1}(V)$ with inverse $t^{-1}\in
ext_{1}^{1}(V),$ i.e., $t^{-1}\circ t=t\circ t^{-1}=i_{V}.$ Then, $%
\underline{(t^{-1})}\in ext(V)$ is the inverse of $\underline{t}\in ext(V),$
i.e., 
\begin{equation}
(\underline{t})^{-1}=\underline{(t^{-1})}.  \label{EO.5}
\end{equation}%
Indeed, by using Eq.(\ref{EO.4}) and the obvious property $\underline{i}%
_{V}=i_{\bigwedge V},$ we have that 
\begin{equation*}
t^{-1}\circ t=t\circ t^{-1}=i_{V}\Rightarrow \underline{(t^{-1})}\circ 
\underline{t}=\underline{t}\circ \underline{(t^{-1})}=i_{\bigwedge V},
\end{equation*}%
which means that the inverse of the extended of $t$ equals the extended of
the inverse of $t.$

In accordance with the above corollary we use in what follows a more simple
notation $\underline{t}^{-1}$ to denote both $(\underline{t})^{-1}$ and $%
\underline{(t^{-1})}.$

Let $\{e_{j}\}$ be any basis for $V,$ and $\{e^{j}\}$ its euclidean
reciprocal basis for $V,$ i.e., $e_{j}\cdot e^{k}=\delta _{j}^{k}.$ There
are two interesting and useful formulas for calculating the extended of $%
t\in ext_{1}^{1}(V),$ i.e.,

\begin{eqnarray}
\underline{t}(X) &=&1\cdot X+\overset{n}{\underset{k=1}{\sum }}\frac{1}{k!}%
(e^{j_{1}}\wedge \ldots \wedge e^{j_{k}})\cdot Xt(e_{j_{1}})\wedge \ldots
\wedge t(e_{j_{k}})  \label{EO.6a} \\
&=&1\cdot X+\overset{n}{\underset{k=1}{\sum }}\frac{1}{k!}(e_{j_{1}}\wedge
\ldots \wedge e_{j_{k}})\cdot Xt(e^{j_{1}})\wedge \ldots \wedge t(e^{j_{k}}).
\label{EO.6b}
\end{eqnarray}

\section{Adjoint Operator}

\subsection{Standard Adjoint Operator}

Let $\bigwedge_{1}^{\diamond }V$ and $\bigwedge_{2}^{\diamond }V$ be two
subspaces of $\bigwedge V$ such that each of them is any sum of homogeneous
subspaces of $\bigwedge V.$ Let $\{e_{j}\}$ and $\{e^{j}\}$ be two euclidean
reciprocal bases to each other for $V,$ i.e., $e_{j}\cdot e^{k}=\delta
_{j}^{k}.$

We call \emph{standard adjoint operator }of $t$ the linear mapping $1$-$%
ext(\bigwedge_{1}^{\diamond }V;\bigwedge_{2}^{\diamond }V)\ni t\rightarrow
t^{\dagger }\in 1$-$ext(\bigwedge_{2}^{\diamond }V;\bigwedge_{1}^{\diamond
}V)$ such that for any $Y\in \bigwedge_{2}^{\diamond }V:$%
\begin{eqnarray}
t^{\dagger }(Y) &=&t(\left\langle 1\right\rangle _{\bigwedge_{1}^{\diamond
}V})\cdot Y+\underset{k=1}{\overset{n}{\sum }}\frac{1}{k!}t(\left\langle
e^{j_{1}}\wedge \ldots e^{j_{k}}\right\rangle _{\bigwedge_{1}^{\diamond
}V})\cdot Ye_{j_{1}}\wedge \ldots e_{j_{k}}  \label{SAO.1a} \\
&=&t(\left\langle 1\right\rangle _{\bigwedge_{1}^{\diamond }V})\cdot Y+%
\underset{k=1}{\overset{n}{\sum }}\frac{1}{k!}t(\left\langle e_{j_{1}}\wedge
\ldots e_{j_{k}}\right\rangle _{\bigwedge_{1}^{\diamond }V})\cdot
Ye^{j_{1}}\wedge \ldots e^{j_{k}}.  \label{SAO.1b}
\end{eqnarray}
Using a more compact notation by employing the \emph{collective index} $J$
we can write 
\begin{eqnarray}
t^{\dagger }(Y) &=&\underset{J}{\sum }\frac{1}{\nu (J)!}t(\left\langle
e^{J}\right\rangle _{\bigwedge_{1}^{\diamond }V})\cdot Ye_{J}  \label{SAO.2a}
\\
&=&\underset{J}{\sum }\frac{1}{\nu (J)!}t(\left\langle e_{J}\right\rangle
_{\bigwedge_{1}^{\diamond }V})\cdot Ye^{J},  \label{SAO.2b}
\end{eqnarray}
We call $t^{\dagger }$ the \emph{standard adjoint} of $t.$ It should be
noticed the use of the $\bigwedge_{1}^{\diamond }V$-projector extensor.

The standard adjoint operator is well-defined since the sums appearing in
each one of the above places do not depend on the choice of $\{e_{j}\}.$

Let us take $X\in \bigwedge_{1}^{\diamond }V$ and $Y\in
\bigwedge_{2}^{\diamond }V.$ A straightforward calculation yields 
\begin{eqnarray*}
X\cdot t^{\dagger }(Y) &=&\underset{J}{\sum }\frac{1}{\nu (J)!}%
t(\left\langle e^{J}\right\rangle _{\bigwedge_{1}^{\diamond }V})\cdot
Y(X\cdot e_{J}) \\
&=&t(\underset{J}{\sum }\frac{1}{\nu (J)!}\left\langle (X\cdot
e_{J})e^{J}\right\rangle _{\bigwedge_{1}^{\diamond }V})\cdot Y \\
&=&t(\left\langle X\right\rangle _{\bigwedge_{1}^{\diamond }V})\cdot Y,
\end{eqnarray*}
i.e., 
\begin{equation}
X\cdot t^{\dagger }(Y)=t(X)\cdot Y.  \label{SAO.3}
\end{equation}
It is a generalization of the well-known property which holds for linear
operators.

Let us take $t\in 1$-$ext(\bigwedge_{1}^{\diamond }V;\bigwedge_{2}^{\diamond
}V)$ and $u\in 1$-$ext(\bigwedge_{2}^{\diamond }V;\bigwedge_{3}^{\diamond
}V).$ We can note that $u\circ t\in 1$-$ext(\bigwedge_{1}^{\diamond
}V;\bigwedge_{3}^{\diamond }V)$ and $t^{\dagger }\circ u^{\dagger }\in 1$-$%
ext(\bigwedge_{3}^{\diamond }V;\bigwedge_{1}^{\diamond }V).$ Then, let us
take $X\in \bigwedge_{1}^{\diamond }V$ and $Z\in \bigwedge_{3}^{\diamond }V,$
by using Eq.(\ref{SAO.3}) we have that 
\begin{equation*}
X\cdot (u\circ t)^{\dagger }(Z)=(u\circ t)(X)\cdot Z=t(X)\cdot u^{\dagger
}(Z)=X\cdot (t^{\dagger }\circ u^{\dagger })(Z).
\end{equation*}%
Hence, we get 
\begin{equation}
(u\circ t)^{\dagger }=t^{\dagger }\circ u^{\dagger }.  \label{SAO.4}
\end{equation}

Let us take $t\in 1$-$ext(\bigwedge^{\diamond }V;\bigwedge^{\diamond }V)$
with inverse $t^{-1}\in 1$-$ext(\bigwedge^{\diamond }V;\bigwedge^{\diamond
}V),$ i.e., $t^{-1}\circ t=t\circ t^{-1}=i_{\bigwedge^{\diamond }V},$ where $%
i_{\bigwedge^{\diamond }V}\in 1$-$ext(\bigwedge^{\diamond
}V;\bigwedge^{\diamond }V)$ is the so-called identity function for $%
\bigwedge^{\diamond }V.$ By using Eq.(\ref{SAO.4}) and the obvious property $%
i_{\bigwedge^{\diamond }V}=i_{\bigwedge^{\diamond }V}^{\dagger },$ we have
that 
\begin{equation*}
t^{-1}\circ t=t\circ t^{-1}=i_{\bigwedge^{\diamond }V}\Rightarrow t^{\dagger
}\circ (t^{-1})^{\dagger }=(t^{-1})^{\dagger }\circ t^{\dagger
}=i_{\bigwedge^{\diamond }V},
\end{equation*}%
hence, 
\begin{equation}
(t^{\dagger })^{-1}=(t^{-1})^{\dagger },  \label{SAO.5}
\end{equation}%
i.e., the inverse of the adjoint of $t$ equals the adjoint of the inverse of 
$t$. In accordance with the above corollary it is possible to use a more
simple symbol, say $t^{\ast }$, to denote both of $(t^{\dagger })^{-1}$ and $%
(t^{-1})^{\dagger }.$

Let us take $t\in ext_{1}^{1}(V).$ We note that $\underline{t}\in ext(V)$
and $\underline{(t^{\dagger })}\in ext(V).$ A straightforward calculation by
using Eqs.(\ref{EO.6a}) and (\ref{EO.6b}) yields 
\begin{eqnarray*}
\underline{(t^{\dagger })}(Y) &=&1\cdot Y+\overset{n}{\underset{k=1}{\sum }}%
\frac{1}{k!}(e^{j_{1}}\wedge \ldots e^{j_{k}})\cdot Yt^{\dagger
}(e_{j_{1}})\wedge \ldots t^{\dagger }(e_{j_{k}}) \\
&=&1\cdot Y+ \\
&&\overset{n}{\underset{k=1}{\sum }}\frac{1}{k!}(e^{j_{1}}\wedge \ldots
e^{j_{k}})\cdot Yt^{\dagger }(e_{j_{1}})\cdot e_{p_{1}}e^{p_{1}}\wedge
\ldots t^{\dagger }(e_{j_{k}})\cdot e_{p_{k}}e^{p_{k}} \\
&=&1\cdot Y+\overset{n}{\underset{k=1}{\sum }}\frac{1}{k!}(e_{j_{1}}\cdot
t(e_{p_{1}})e^{j_{1}}\wedge \ldots e_{j_{k}}\cdot
t(e_{p_{k}})e^{j_{k}})\cdot Ye^{p_{1}}\wedge \ldots e^{p_{k}} \\
&=&\underline{t}(1)\cdot Y+\overset{n}{\underset{k=1}{\sum }}\frac{1}{k!}%
\underline{t}(e_{p_{1}}\wedge \ldots e_{p_{k}})\cdot Ye^{p_{1}}\wedge \ldots
e^{p_{k}} \\
&=&(\underline{t})^{\dagger }(Y).
\end{eqnarray*}%
Hence, we get 
\begin{equation}
\underline{(t^{\dagger })}=(\underline{t})^{\dagger }.  \label{SAO.6}
\end{equation}%
This means that the extension operator commutes with the adjoint operator.
In accordance with the above property we may use a more simple notation $%
\underline{t}^{\dagger }$ to denote without ambiguities both $\underline{%
(t^{\dagger })}$ and $(\underline{t})^{\dagger }.$

\subsection{Metric Adjoint Operator}

As we know \cite{2} whenever $V$ is endowed with another metric $G$ (besides 
$G_{E}$) there exists an unique $(1,1)$-extensor $g$ such that the $G$%
-scalar product of $X,Y\in \bigwedge V,$ namely $X\underset{g}{\cdot }Y,$ is
given by 
\begin{equation}
X\underset{g}{\cdot }Y=\underline{g}(X)\cdot Y.  \label{MAO.1}
\end{equation}%
Such $g\in ext_{1}^{1}(V)$ is symmetric and non-degenerate, and has
signature $(p,q),$ i.e., $g=g^{\dagger },$ $\det [g]\neq 0,$ $g$ has $p$
positive and $q$ negative ($p+q=n$) eigenvalues. It is called the metric
extensor for $G.$

To each $t\in 1$-$ext(\bigwedge_{1}^{\diamond }V;\bigwedge_{2}^{\diamond }V)$
we can assign $t^{\dagger (g)}\in 1$-$ext(\bigwedge_{2}^{\diamond
}V;\bigwedge_{1}^{\diamond }V)$ defined as follows 
\begin{equation}
t^{\dagger (g)}=\underline{g}^{-1}\circ t^{\dagger }\circ \underline{g}.
\label{MAU.2}
\end{equation}
It will be called the\emph{\ metric adjoint} of $t.$

As we can easily see, $t^{\dagger (g)}$ is the unique extensor from $%
\bigwedge_{2}^{\diamond }V$ to $\bigwedge_{1}^{\diamond }V$ which satisfies
the fundamental property 
\begin{equation}
X\underset{g}{\cdot }t^{\dagger (g)}(Y)=t(X)\underset{g}{\cdot }Y,
\label{MAO.3}
\end{equation}
for all $X\in \bigwedge_{1}^{\diamond }V$ and $Y\in \bigwedge_{2}^{\diamond
}V.$

The noticeable property given by Eq.(\ref{MAO.3}) is the \emph{metric version%
} of the fundamental property given by Eq.(\ref{SAO.3}).

\section{The Generalization Operator}

\subsection{Standard Generalization Operator}

Let $\{e_{k}\}$ be any basis for $V,$ and $\{e^{k}\}$ be its euclidean
reciprocal basis for $V,$ as we know, $e_{k}\cdot e^{l}=\delta _{k}^{l}.$

The linear mapping $ext_{1}^{1}(V)\ni t\mapsto \underset{\thicksim }{t}\in
ext(V)$ such that for any $X\in \bigwedge V$ 
\begin{equation}
\underset{\thicksim }{t}(X)=t(e^{k})\wedge (e_{k}\lrcorner X)=t(e_{k})\wedge
(e^{k}\lrcorner X)  \label{SGO.1}
\end{equation}
will be called the \emph{generalization operator.} We call $\underset{%
\thicksim }{t}$ the \emph{generalized of }$t.$

The generalization operator is well-defined since it does not depend on the
choice of $\{e_{k}\}.$

We present now some important properties which are satisfied by the
generalization operator.\medskip

\textbf{g1 }The generalization operator is grade-preserving, i.e., 
\begin{equation}
\text{if }X\in \bigwedge^{k}V,\text{ then }\underset{\thicksim }{t}(X)\in
\bigwedge^{k}V.  \label{SGO.1a}
\end{equation}%
\medskip

\textbf{g2 }The grade involution $\widehat{\left. {}\right. }\in ext(V),$
reversion $\widetilde{\left. {}\right. }\in ext(V),$ and conjugation $%
\overline{\left. {}\right. }\in ext(V)$ commute with the generalization
operator, i.e., 
\begin{eqnarray}
\underset{\thicksim }{t}(\widehat{X}) &=&\widehat{\underset{\thicksim }{t}(X)%
},  \label{SGO.2a} \\
\underset{\thicksim }{t}(\widetilde{X}) &=&\widetilde{\underset{\thicksim }{t%
}(X)},  \label{SGO.2b} \\
\underset{\thicksim }{t}(\overline{X}) &=&\overline{\underset{\thicksim }{t}%
(X)}.  \label{SGO.2c}
\end{eqnarray}%
They are immediate consequences of the grade-preserving property.\medskip

\textbf{g3 }For any $\alpha \in \mathbb{R},$ $v\in V$ and $X,Y\in \bigwedge
V $ it holds 
\begin{eqnarray}
\underset{\thicksim }{t}(\alpha ) &=&0,  \label{SGO.3a} \\
\underset{\thicksim }{t}(v) &=&t(v),  \label{SGO.3b} \\
\underset{\thicksim }{t}(X\wedge Y) &=&\underset{\thicksim }{t}(X)\wedge
Y+X\wedge \underset{\thicksim }{t}(Y).  \label{SGO.3c}
\end{eqnarray}

The proof of Eq.(\ref{SGO.3a}) and Eq.(\ref{SGO.3b}) are left to the reader.
Hint: $v\lrcorner \alpha =0$ and $v\lrcorner w=v\cdot w.$ Now, the
identities: $a\lrcorner (X\wedge Y)=(a\lrcorner X)\wedge Y+\widehat{X}\wedge
(a\lrcorner Y)$ and $a\wedge X=\widehat{X}\wedge a,$ with $a\in V$ and $%
X,Y\in \bigwedge V,$ allow us to prove the property given by Eq.(\ref{SGO.3c}%
).

We can prove that the basic properties given by Eq.(\ref{SGO.3a}), Eq.(\ref%
{SGO.3b}) and Eq.(\ref{SGO.3c}) together are completely equivalent to the 
\emph{generalization procedure} as defined by Eq.(\ref{SGO.1}).\medskip

\textbf{g4} The generalization operator commutes with the adjoint operator,
i.e., 
\begin{equation}
(\underset{\thicksim }{t})^{\dagger }=\underset{\thicksim }{(t^{\dagger })},
\label{SGO.4}
\end{equation}%
or put it on another way, the adjoint of the generalized of $t$ is just the
generalized of the adjoint of $t$.\medskip

\textbf{Proof}

A straightforward calculation by using Eq.(\ref{SAO.3}) and the multivector
identities: $X\cdot (a\wedge Y)=(a\lrcorner X)\wedge Y$ and $X\cdot
(a\lrcorner Y)=(a\wedge X)\cdot Y,$ with $a\in V$ and $X,Y\in \bigwedge V,$
gives 
\begin{eqnarray*}
(\underset{\thicksim }{t})^{\dagger }(X)\cdot Y &=&X\cdot \underset{%
\thicksim }{t}(Y) \\
&=&(e_{j}\wedge (t(e^{j})\lrcorner X))\cdot Y=(e_{j}\wedge (t(e^{j})\cdot
e^{k}e_{k}\lrcorner X))\cdot Y \\
&=&(e^{j}\cdot t^{\dagger }(e^{k})e_{j}\wedge (e_{k}\lrcorner X))\cdot
Y=(t^{\dagger }(e^{k})\wedge (e_{k}\lrcorner X))\cdot Y \\
&=&\underset{\thicksim }{(t^{\dagger })}(X)\cdot Y.
\end{eqnarray*}%
Hence, by the non-degeneracy property of the euclidean scalar product, the
required result follows.$\blacksquare $

In agreement with the above property we use in what follows a more simple
symbol, $\underset{\thicksim }{t^{\dagger }}$ to denote both $(\underset{%
\thicksim }{t})^{\dagger }$ or $\underset{\thicksim }{(t^{\dagger })}.$

\textbf{g5} The symmetric (skew-symmetric) part of the generalized of $t$ is
just the generalized of the symmetric (skew-symmetric) part of $t,$ i.e., 
\begin{equation}
(\underset{\thicksim }{t})_{\pm }=(\underset{\thicksim }{t_{\pm }}).
\label{SGO.5}
\end{equation}%
This property follows immediately from Eq.(\ref{SGO.4}).

We see also that it is possible to use a more simple notation, $\underset{%
\thicksim }{t}_{\pm }$to denote $(\underset{\thicksim }{t})_{\pm }$ or $(%
\underset{\thicksim }{t_{\pm }}).$

\textbf{g6} The skew-symmetric part of the generalized of $t$ can be
factorized by the noticeable formula\footnote{%
Recall that $X\times Y\equiv \frac{1}{2}(XY-YX).$} 
\begin{equation}
\underset{\thicksim }{t}_{-}(X)=\frac{1}{2}biv[t]\times X,  \label{SGO.6}
\end{equation}%
where $biv[t]\equiv t(e^{k})\wedge e_{k}$ is a \emph{characteristic invariant%
} of $t$, called the \emph{bivector of} $t.\medskip $

\textbf{Proof}

By using Eq.(\ref{SGO.5}), the well-known identity $t_{-}(a)=\frac{1}{2}%
biv[t]\times a$ and the multivector identity $B\times X=(B\times
e^{k})\wedge (e_{k}\lrcorner X),$ with $B\in \bigwedge^{2}V$ and $X\in
\bigwedge V,$ we have that 
\begin{equation*}
\underset{\thicksim }{t}_{-}(X)=t_{-}(e^{k})\wedge (e_{k}\lrcorner X)=(\frac{%
1}{2}biv[t]\times e^{k})\wedge (e_{k}\lrcorner X)=\frac{1}{2}biv[t]\times
X.\blacksquare
\end{equation*}

\textbf{g7} A noticeable formula holds for the skew-symmetric part of the
generalized of $t.$ For all $X,Y\in \bigwedge V$ 
\begin{equation}
\underset{\thicksim }{t}_{-}(X\ast Y)=\underset{\thicksim }{t}_{-}(X)\ast
Y+X\ast \underset{\thicksim }{t}_{-}(Y),  \label{SGO.7}
\end{equation}%
where $\ast $ is any product either $(\wedge ),$ $(\cdot ),$ $(\lrcorner
,\llcorner )$ or $($\emph{Clifford product}$).$

In order to prove this property we must use Eq.(\ref{SGO.6}) and the
multivector identity $B\times (X\ast Y)=(B\times X)\ast Y+X\ast (B\times Y),$
with $B\in \bigwedge^{2}V$ and $X,Y\in \bigwedge V$. By taking into account
Eq.(\ref{SGO.3a}) we can see that the following property for the euclidean
scalar product of multivectors holds 
\begin{equation}
\underset{\thicksim }{t}_{-}(X)\cdot Y+X\cdot \underset{\thicksim }{t}%
_{-}(Y)=0.  \label{SGO.7a}
\end{equation}%
It is consistent with the well-known property: \emph{the adjoint of a
skew-symmetric extensor equals minus the extensor!}

\section{Determinant}

Let $t$ be any $(1,1)$-extensor. We define the \emph{determinant}\footnote{%
The concept of determinant of a $(1,1)$-extensor is related, but distinct
from the well known determinant of a square matrix. For details the reader
is invited to consult \cite{rodoliv2006}.}\emph{\ of }$t$ as the unique real
number, denoted by $\det [t],$ such that for all non-zero pseudoscalar $I$ 
\begin{equation}
\underline{t}(I)=\det [t]I.  \label{D.1}
\end{equation}

It is a \emph{well-defined }scalar \emph{invariant} since it does not depend
on the choice of $I.$

We present now some of the most important properties satisfied by the
determinant.

\textbf{d1} Let $t$ and $u$ be two $(1,1)$--extensors. It holds 
\begin{equation}
\det [u\circ t]=\det [u]\det [t].  \label{D.2}
\end{equation}

\textbf{Proof}

Take a non-zero pseudoscalar $I\in \bigwedge^{n}V.$ By using Eq.(\ref{EO.4})
and Eq.(\ref{D.1}) we can write that 
\begin{eqnarray*}
\det [u\circ t]I &=&\underline{u\circ t}(I)=\underline{u}\circ \underline{t}%
(I)=\underline{u}(\underline{t}(I)) \\
&=&\underline{u}(\det [t]I)=\det [t]\underline{u}(I), \\
&=&\det [t]\det [u]I.\blacksquare
\end{eqnarray*}

\textbf{d2} Let us take $t\in ext_{1}^{1}(V)$ with inverse $t^{-1}\in
ext_{1}^{1}(V).$ It holds 
\begin{equation}
\det [t^{-1}]=(\det [t])^{-1}.  \label{D.3}
\end{equation}

Indeed, by using Eq.(\ref{D.2}) and the obvious property $\det [i_{V}]=1,$
we have that 
\begin{equation*}
t^{-1}\circ t=t\circ t^{-1}=i_{V}\Rightarrow \det [t^{-1}]\det [t]=\det
[t]\det [t^{-1}]=1,
\end{equation*}%
which means that the \emph{determinant of the inverse equals the inverse of
the determinant.}

Due to the above corollary it is often convenient to use the short notation $%
\left. \det \right. ^{-1}[t]$ for both $\det [t^{-1}]$ and $(\det [t])^{-1}. 
$

\textbf{d3} Let us take $t\in ext_{1}^{1}(V).$ It holds 
\begin{equation}
\det [t^{\dagger }]=\det [t].  \label{D.4}
\end{equation}%
Indeed, take a non-zero $I\in \bigwedge^{n}V.$ Then, by using Eq.(\ref{D.1})
and Eq.(\ref{SAO.3}) we have that 
\begin{equation*}
\det [t^{\dagger }]I\cdot I=\underline{t}^{\dagger }(I)\cdot I=I\cdot 
\underline{t}(I)=I\cdot \det [t]I=\det [t]I\cdot I,
\end{equation*}%
whence, the expected result follows.

Let $\{e_{j}\}$ be any basis for $V,$ and $\{e^{j}\}$ be its euclidean
reciprocal basis for $V,$ i.e., $e_{j}\cdot e^{k}=\delta _{j}^{k}.$ There
are two interesting and useful formulas for calculating $\det [t],$ i.e., 
\begin{eqnarray}
\det [t] &=&\underline{t}(e_{1}\wedge \ldots \wedge e_{n})\cdot (e^{1}\wedge
\ldots \wedge e^{n}),  \label{D.5a} \\
&=&\underline{t}(e^{1}\wedge \ldots \wedge e^{n})\cdot (e_{1}\wedge \ldots
\wedge e_{n}).  \label{D.5b}
\end{eqnarray}%
They follow from Eq.(\ref{D.1}) by using $(e_{1}\wedge \ldots \wedge
e_{n})\cdot (e^{1}\wedge \ldots \wedge e^{n})=1$ which is an immediate
consequence of the formula for the euclidean scalar product of simple $k$%
-vectors and the reciprocity property of $\{e_{k}\}$ and $\{e^{k}\}$.

Each of Eq.(\ref{D.5a}) and Eq.(\ref{D.5b}) is completely equivalent to the
definition of determinant given by Eq.(\ref{D.1}).

We will end this section presenting an useful formula for the inversion of a
non-singular $(1,1)$--extensor.

Let us take $t\in ext_{1}^{1}(V)$. If $t$ is non-singular, i.e., $\det
[t]\neq 0,$ then there exists its inverse $t^{-1}\in ext_{1}^{1}(V)$ which
is given by 
\begin{equation}
t^{-1}(v)=\left. \det \right. ^{-1}[t]\underline{t}^{\dagger }(vI)I^{-1},
\label{D.6}
\end{equation}
where $I\in \bigwedge^{n}V$ is any non-zero pseudoscalar.

\textbf{Proof}

We must prove that $t^{-1}$ given by the formula above satisfies both of
conditions $t^{-1}\circ t=i_{V}$ and $t\circ t^{-1}=i_{V}.$

Let $I\in \bigwedge^{n}V$ be a non-zero pseudoscalar. Take $v\in V,$ by
using the extensor identities\footnote{%
These extensor identities follow directly from the fundamental identity $%
X\lrcorner \underline{t}(Y)=\underline{t}(\underline{t}^{\dagger
}(X)\lrcorner Y)$ with $X,Y\in \bigwedge V$. For the first one: take $X=v,$ $%
Y=I$ and use $(t^{\dagger })^{\dagger }=t,$ eq.(\ref{D.1}) and $\det
[t^{\dagger }]=\det [t].$ For the second one: take $X=vI,Y=I^{-1}$ and use
eq.(\ref{D.1}).} $\underline{t}^{\dagger }(t(v)I)I^{-1}=t(\underline{t}%
^{\dagger }(vI)I^{-1})=\det [t]v,$ we have that 
\begin{equation*}
t^{-1}\circ t(v)=t^{-1}(t(v))=\left. \det \right. ^{-1}[t]\underline{t}%
^{\dagger }(t(v)I)I^{-1}=\left. \det \right. ^{-1}[t]\det [t]v=i_{V}(v).
\end{equation*}%
And 
\begin{equation*}
t\circ t^{-1}(v)=t(t^{-1}(v))=\left. \det \right. ^{-1}[t]t(\underline{t}%
^{\dagger }(vI)I^{-1})=\left. \det \right. ^{-1}[t]\det
[t]v=i_{V}(v).\blacksquare
\end{equation*}

\section{Hodge Extensor}

\subsection{Standard Hodge Extensor}

Let $\{e_{j}\}$ and $\{e^{j}\}$ be two euclidean reciprocal bases to each
other for $V,$ i.e., $e_{j}\cdot e^{k}=\delta _{j}^{k}.$ Associated to them
we define a non-zero pseudoscalar 
\begin{equation}
\tau =\sqrt{e_{\wedge }\cdot e_{\wedge }}e^{\wedge },  \label{SHE.1}
\end{equation}
where $e_{\wedge }\equiv e_{1}\wedge \ldots \wedge e_{n}\in \bigwedge^{n}V$
and $e^{\wedge }\equiv e^{1}\wedge \ldots \wedge e^{n}\in \bigwedge^{n}V$.
Note that $e_{\wedge }\cdot e_{\wedge }>0,$ since an euclidean scalar
product is positive definite. Such $\tau $ will be called a \emph{standard
volume pseudoscalar} for $V.$

The standard volume pseudoscalar has the fundamental property 
\begin{equation}
\tau \cdot \tau =\tau \lrcorner \widetilde{\tau }=\tau \widetilde{\tau }=1,
\label{SHE.2}
\end{equation}
which follows from the obvious result $e_{\wedge }\cdot e^{\wedge }=1$.

From Eq.(\ref{SHE.2}), we can easily get an expansion formula for
pseudoscalars of $\bigwedge^{n}V,$ i.e., 
\begin{equation}
I=(I\cdot \tau )\tau .  \label{SHE.3}
\end{equation}

The extensor $\star \in ext(V)$ which is defined by $\star :\bigwedge
V\rightarrow \bigwedge V$ such that 
\begin{equation}
\star X=\widetilde{X}\lrcorner \tau =\widetilde{X}\tau ,  \label{SHE.4}
\end{equation}
will be called a \emph{standard Hodge extensor }on $V.$

It should be noticed that 
\begin{equation}
\text{if }X\in \bigwedge^{p}V,\text{ then }\star X\in \bigwedge^{n-p}V.
\label{SHE.4a}
\end{equation}
That means that $\star $ can be also defined as a $(p,n-p)$-extensor over $%
V. $

The extensor over $V,$ namely $\star ^{-1},$ which is given by $\star
^{-1}:\bigwedge V\rightarrow \bigwedge V$ such that 
\begin{equation}
\star ^{-1}X=\tau \llcorner \widetilde{X}=\tau \widetilde{X}  \label{SHE.5}
\end{equation}
is the \emph{inverse extensor} of $\star .$

Let us take $X\in \bigwedge V$. By Eq.(\ref{SHE.2}), we have indeed that $%
\star ^{-1}\circ \star X=\tau \widetilde{\tau }X=X,$ and $\star \circ \star
^{-1}X=X\widetilde{\tau }\tau =X,$ i.e., $\star ^{-1}\circ \star =\star
\circ \star ^{-1}=i_{\bigwedge V},$ where $i_{\bigwedge V}\in ext(V)$ is the
so-called \emph{identity function }for $\bigwedge V.$

Let us take $X,Y\in \bigwedge V.$ By using the multivector identity $%
(XA)\cdot Y=X\cdot (Y\widetilde{A})$ and Eq.(\ref{SHE.2}) we get 
\begin{equation}
(\star X)\cdot (\star Y)=X\cdot Y.  \label{SHE.6}
\end{equation}%
That means that the standard Hodge extensor preserves the euclidean scalar
product.

Let us take $X,Y\in \bigwedge^{p}V.$ By using Eq.(\ref{SHE.3}) together with
the multivector identity $(X\wedge Y)\cdot Z=Y\cdot (\widetilde{X}\lrcorner
Z),$ and Eq.(\ref{SHE.6}) we get 
\begin{equation}
X\wedge (\star Y)=(X\cdot Y)\tau .  \label{SHE.7}
\end{equation}%
This noticeable identity is completely equivalent to the definition of the
standard Hodge extensor given by Eq.(\ref{SHE.4}).

Let us take $X\in \bigwedge^{p}V$ and $Y\in \bigwedge^{n-p}V.$ By using the
multivector identity $(X\lrcorner Y)\cdot Z=Y\cdot (\widetilde{X}\wedge Z)$
and Eq.(\ref{SHE.3}) we get 
\begin{equation}
(\star X)\cdot Y\tau =X\wedge Y.  \label{SHE.8}
\end{equation}

\subsection{Metric Hodge Extensor}

Let $g$ be a metric extensor on $V$ with signature $(p,q),$ i.e., $g\in
ext_{1}^{1}(V)$ such that $g=g^{\dagger }$ and $\det [g]\neq 0,$ and it has $%
p$ positive and $q$ negative eigenvalues. Associated to $\{e_{j}\}$ and $%
\{e^{j}\}$ we can define another non-zero pseudoscalar 
\begin{equation}
\underset{g}{\tau }=\sqrt{\left\vert e_{\wedge }\underset{g}{\cdot }%
e^{\wedge }\right\vert }e^{\wedge }=\sqrt{\left\vert \det [g]\right\vert }%
\tau .  \label{MHE.1}
\end{equation}%
It will be called a metric volume pseudoscalar for $V.$ It has the
fundamental property 
\begin{equation}
\underset{g}{\tau }\underset{g^{-1}}{\cdot }\underset{g}{\tau }=\underset{g}{%
\tau }\underset{g^{-1}}{\lrcorner }\widetilde{\underset{g}{\tau }}=\underset{%
g}{\tau }\underset{g^{-1}}{}\widetilde{\underset{g}{\tau }}=(-1)^{q}.
\label{MHE.2}
\end{equation}%
Eq.(\ref{MHE.2}) follows from Eq.(\ref{SHE.2}) by taking into account the
definition of determinant of a $(1,1)$-extensor, and recalling that $%
sgn(\det [g])=(-1)^{q}.$

An expansion formula for pseudoscalars of $\bigwedge^{n}V$ can be also
obtained from Eq.(\ref{MHE.2}), i.e., 
\begin{equation}
I=(-1)^{q}(I\underset{g^{-1}}{\cdot }\underset{g}{\tau })\underset{g}{\tau }.
\label{MHE.3}
\end{equation}

The extensor $\underset{g}{\star }\in ext(V)$ which is defined by $\underset{%
g}{\star }:\bigwedge V\rightarrow \bigwedge V$ such that 
\begin{equation}
\underset{g}{\star }X=\widetilde{X}\underset{g^{-1}}{\lrcorner }\underset{g}{%
\tau }=\widetilde{X}\underset{g^{-1}}{}\underset{g}{\tau }  \label{MHE.4}
\end{equation}
will be called a \emph{metric Hodge extensor }on $V.$ It should be noticed
that in general we need to use of both the $g$ and $g^{-1}$ metric Clifford
algebras.

We see that 
\begin{equation}
\text{if }X\in \bigwedge^{p}V,\text{ then }\underset{g}{\star }X\in
\bigwedge^{n-p}V.  \label{MHE.4a }
\end{equation}
It means that $\underset{g}{\star }\in ext(V)$ can also be defined as $%
\underset{g}{\star }\in ext_{p}^{n-p}(V).$

The extensor over $V,$ namely $\underset{g}{\star }^{-1},$ which is given by 
$\underset{g}{\star }^{-1}:\bigwedge V\rightarrow \bigwedge V$ such that 
\begin{equation}
\underset{g}{\star }^{-1}X=(-1)^{q}\underset{g}{\tau }\underset{g^{-1}}{%
\llcorner }\widetilde{X}=(-1)^{q}\underset{g}{\tau }\underset{g^{-1}}{}%
\widetilde{X}  \label{MHE.5}
\end{equation}
is the \emph{inverse extensor} of $\underset{g}{\star }.$

Let us take $X\in \bigwedge V.$ By using Eq.(\ref{MHE.2}), we have indeed
that $\underset{g}{\star }^{-1}\circ \underset{g}{\star }X=(-1)^{q}\underset{%
g}{\tau }\underset{g^{-1}}{}\widetilde{\underset{g}{\tau }}\underset{g^{-1}}{%
}X=X,$ and $\underset{g}{\star }\circ \underset{g}{\star }^{-1}X=(-1)^{q}X%
\underset{g^{-1}}{}\widetilde{\underset{g}{\tau }}\underset{g^{-1}}{}%
\underset{g}{\tau }=X,$ i.e., $\underset{g}{\star }^{-1}\circ \star =%
\underset{g}{\star }\circ \underset{g}{\star }^{-1}=i_{\bigwedge V}.$

Take $X,Y\in \bigwedge V.$ The identity $(X\underset{g^{-1}}{}A)\underset{%
g^{-1}}{\cdot }Y=X\underset{g^{-1}}{\cdot }(Y\underset{g^{-1}}{}\widetilde{A}%
)$ and Eq.(\ref{MHE.2}) yield 
\begin{equation}
(\underset{g}{\star }X)\underset{g^{-1}}{\cdot }(\underset{g}{\star }%
Y)=(-1)^{q}X\underset{g^{-1}}{\cdot }Y.  \label{MHE.6}
\end{equation}

Take $X,Y\in \bigwedge^{p}V.$ Eq.(\ref{MHE.3}), the identity $(X\wedge Y)%
\underset{g^{-1}}{\cdot }Z=Y\underset{g^{-1}}{\cdot }(\widetilde{X}\underset{%
g^{-1}}{\lrcorner }Z)$ and Eq.(\ref{MHE.6}) allow us to obtain 
\begin{equation}
X\wedge (\underset{g}{\star }Y)=(X\underset{g^{-1}}{\cdot }Y)\underset{g}{%
\tau }.  \label{MHE.7}
\end{equation}%
This remarkable property is completely equivalent to the definition of the
metric Hodge extensor given by Eq.(\ref{MHE.4}).

Take $X\in \bigwedge^{p}V$ and $Y\in \bigwedge^{n-p}V.$ The use of identity $%
(X\underset{g^{-1}}{\lrcorner }Y)\underset{g^{-1}}{\cdot }Z=Y\underset{g^{-1}%
}{\cdot }(\widetilde{X}\wedge Z)$ and Eq.(\ref{MHE.3}) yield 
\begin{equation}
(\underset{g}{\star }X)\underset{g^{-1}}{\cdot }Y\underset{g}{\tau }%
=(-1)^{q}X\wedge Y.  \label{MHE.8}
\end{equation}

It might as well be asked what is the relationship between the standard and
metric Hodge extensors as defined above by Eq.(\ref{SHE.4}) and Eq.(\ref%
{MHE.4}).

Take $X\in \bigwedge V.$ By using Eq.(\ref{MHE.1}), the multivector identity
for an invertible $(1,1)$-extensor $\underline{t}^{-1}(X)\lrcorner Y=%
\underline{t}^{\dagger }(X\lrcorner \underline{t}^{\ast }(Y)),$ and the
definition of determinant of a $(1,1)$-extensor we have that 
\begin{eqnarray*}
\underset{g}{\star }X &=&\underline{g}^{-1}(\widetilde{X})\lrcorner \sqrt{%
\left\vert \det [g]\right\vert }\tau =\sqrt{\left\vert \det [g]\right\vert }%
\underline{g}(\widetilde{X}\lrcorner \underline{g}^{-1}(\tau )) \\
&=&\frac{\sqrt{\left\vert \det [g]\right\vert }}{\det [g]}\underline{g}(%
\widetilde{X}\lrcorner \tau )=\frac{sgn(\det [g])}{\sqrt{\left\vert \det
[g]\right\vert }}\underline{g}\circ \star (X),
\end{eqnarray*}%
i.e., 
\begin{equation}
\underset{g}{\star }=\frac{(-1)^{q}}{\sqrt{\left\vert \det [g]\right\vert }}%
\underline{g}\circ \star .  \label{MHE.9}
\end{equation}%
Eq.(\ref{MHE.9}) is then the formula which relates the metric Hodge extensor 
$\underset{g}{\star }$ with the standard Hodge extensor $\star .$

We know \cite{2} that for any metric extensor $g\in ext_{1}^{1}(V)$ there
exists a non-singular $(1,1)$-extensor $h\in ext_{1}^{1}(V)$ such that 
\begin{equation}
g=h^{\dagger }\circ \eta \circ h,  \label{MHE.10}
\end{equation}%
where $\eta \in ext_{1}^{1}(V)$ is an orthogonal metric extensor with the
same signature as $g$. Such $h$ is called a \emph{gauge extensor} for $g.$

We can also get a noticeable formula which relates the $g$-metric Hodge
extensor with the $\eta $-metric Hodge extensor.

As we know, the $g$ and $g^{-1}$ contracted products $\underset{g}{\lrcorner 
}$ and $\underset{g^{-1}}{\lrcorner }$ are related to the $\eta $-contracted
product$\underset{\eta }{\lrcorner }$ (recall that $\eta =\eta ^{-1}$) by
the following \emph{golden formulas} 
\begin{eqnarray}
\underline{h}(X\underset{g}{\lrcorner }Y) &=&\underline{h}(X)\underset{\eta }%
{\lrcorner }\underline{h}(Y),  \label{MHE.10a} \\
\underline{h}^{*}(X\underset{g^{-1}}{\lrcorner }Y) &=&\underline{h}^{*}(X)%
\underset{\eta }{\lrcorner }\underline{h}^{*}(Y).  \label{MHE.10b}
\end{eqnarray}

Now, take $X\in \bigwedge V.$ By using Eq.(\ref{MHE.10b}), Eq.(\ref{MHE.1}),
the definition of determinant of a $(1,1)$-extensor, Eq.(\ref{MHE.10}) and
the obvious equation $\underset{\eta }{\tau }=\tau $ we have that 
\begin{eqnarray*}
\underset{g}{\star }X &=&\underline{h}^{\dagger }(\underline{h}^{\ast }(%
\widetilde{X})\underset{\eta }{\lrcorner }\underline{h}^{\ast }(\underset{g}{%
\tau }))=\sqrt{\left\vert \det [g]\right\vert }\underline{h}^{\dagger }(%
\underline{h}^{\ast }(\widetilde{X})\underset{\eta }{\lrcorner }\det
[h^{\ast }]\tau ) \\
&=&\left\vert \det [h]\right\vert \det [h^{\ast }]\underline{h}^{\dagger }(%
\widetilde{\underline{h}^{\ast }(X)}\underset{\eta ^{-1}}{\lrcorner }%
\underset{\eta }{\tau })=sgn(\det [h])\underline{h}^{\dagger }\circ \underset%
{\eta }{\star }\circ \underline{h}^{\ast }(X),
\end{eqnarray*}%
i.e., 
\begin{equation}
\underset{g}{\star }\text{ }=sgn(\det [h])\underline{h}^{\dagger }\circ 
\underset{\eta }{\star }\circ \underline{h}^{\ast }.  \label{MHE.11}
\end{equation}%
This formula which relates the $g$-metric Hodge extensor $\underset{g}{\star 
}$ with the $\eta $-metric Hodge extensor $\underset{\eta }{\star }$ will
play an important role in the applications we have in mind.

\section{Conclusions}

In this third paper in a series of eight we recalled some basic notions of
the theory of extensors. Together with \cite{1,2} it completes the algebraic
part of our enterprise. The $k$-extensors and $(p,q)$-extensors are
introduced in Section 2 and projector operators are studied in Section 3 .
The extension operator (or exterior power operator) of a given linear
operator $t$ (i.e., a $(1,1)$-extensor) on a real vector space $V$ has been
studied in details in Section 4. In Section 5 the standard and metric
adjoint operators have been given and Section 6 has been dedicated to the
generalization operator of $t$ and its properties. Such generalization
operator plays an important role in our formulation of the differential
geometry of arbitrary (smooth) manifolds which is presented in sequel papers
of this series. Section 8 showed some applications of the concept of
extensors. In particular, the standard and metric Hodge operators are
introduced and the non trivial relation between them is disclosed (Eq.(\ref%
{MHE.9})). A formula (Eq.(\ref{MHE.11})) relating metric Hodge operators
corresponding to deformed metric structures is also derived. These formulas,
every student of General Relativity and modern geometric theories of Physics
will regonize as very useful ones, forthey simplify many involved
calculations (see, e.g., \cite{3,4,5,6,7}). To end, we recall that the
theory of extensors is a part of a more genral theory of multivector
functions and functionals, which are prsented in \cite{8,9}\medskip

\textbf{Acknowledgments: }V. V. Fern\'{a}ndez and A. M. Moya are very
grateful to Mrs. Rosa I. Fern\'{a}ndez who gave to them material and
spiritual support at the starting time of their research work. This paper
could not have been written without her inestimable help.

\end{document}